\documentclass[12pt]{amsart}
\usepackage{amssymb,enumerate,mathrsfs, amsmath, amsthm}
\usepackage{url}

\newtheorem{theorem}{Theorem}[section]
\newtheorem{cor}[theorem]{Corollary}
\newtheorem{lemma}[theorem]{Lemma}

\newtheorem{prop}[theorem]{Proposition}

\newtheorem{example}[theorem]{Example}

\def\beq{\begin{equation}}
\def\eeq{\end{equation}}

\def\tm{\tilde{m}}
\def\tp{\tilde{p}}

\def\tN{\widetilde{N}}
\def\tM{\widetilde{M}}

\numberwithin{equation}{section}

\def\beq{\begin{equation}}
\def\eeq{\end{equation}}

\newcommand{\monsub}[3]{{#1} |_{#2 \rightarrow #3}}
\newcommand{\monsubtwo}[5]{{#1} |_{#2 \rightarrow #3, #4 \rightarrow #5}}

\def\cH{ {{\mathcal H}}}

\def\cP{ {{\mathcal P}}}

\def\bbR{ {\mathbb R}}

\def\gtupn{(\mathbb{R}^{n\times n})^g}

\def\cH{ {\mathcal H} }

\def\cP{{\mathcal P}}

\def\cR{{ \mathcal R }}

\def\rank{\mbox{rank}}

\def\beq{\begin{equation}}
\def\eeq{\end{equation}}

\def\tf{{\widetilde f}}

\def\tM{\widetilde M}
\def\span{\mathrm{span}}

\def\cR{{\mathcal R}}

\newcommand{\purehered}[1]{\cH_{hered}(#1)}
\newcommand{\pureantihered}[1]{\cH_{antihered}(#1)}

\def\tf{\tilde{f}}
\def\tk{\tilde{k}}

\title[Noncommutative Plurisubharmonic Polynomials: Part I]
{Noncommutative Plurisubharmonic Polynomials \\ Part I: Global Assumptions}

\author[Greene]{Jeremy M. Greene${}^1$${}^\dagger$}
\address{Jeremy M. Greene, Department of Mathematics \\
University of California \\
San Diego
}
\email{j1greene@math.ucsd.edu}
\thanks{${}^1$Research supported by NSF grants
DMS-0700758, DMS-0757212.}

\author[Helton]{J. William Helton${}^2$}
\address{J. William Helton, Department of Mathematics\\
  University of California \\
  San Diego}
\email{helton@math.ucsd.edu}
\thanks{${}^2$Research supported by NSF grants
DMS-0700758, DMS-0757212, and the Ford Motor Co.}

\author[Vinnikov]{Victor Vinnikov${}^3$}
\address{Victor Vinnikov, Department of Mathematics\\
  Ben Gurion University of Negev \\ 
  Beer Sheva 84105, Israel}
\email{vinnikov@cs.bgu.ac.il}
\thanks{${}^3$Partially supported by a grant from the Israel Science Foundation.}

\thanks{${}^\dagger$The material in this paper is part of the Ph.D. thesis of 
Jeremy M. Greene at UCSD}


  \makeatletter
    \newcommand{\mycontentsbox}{%
        {\centerline{NOT FOR PUBLICATION}
                  \tableofcontents}}
                    \def\enddoc@text{\ifx\@empty\@translators \else\@settranslators\fi
                        \ifx\@empty\addresses \else\@setaddresses\fi
                            \newpage\mycontentsbox}
                              \makeatother

\begin{document}
\maketitle

\begin{abstract}
We consider symmetric polynomials, $p$, in the noncommutative (nc) free variables 
$\{x_1, x_2, \ldots, x_g\}$. We define the nc complex hessian of $p$ as the 
second directional derivative (replacing $x^T$ by $y$)
$$
q(x,x^T)[h,h^T] := \frac{\partial^2 p}{\partial s \partial t} (x+th, y+sk) |_{t,s=0} |_{y=x^T, k=h^T}.
$$
We call an nc symmetric polynomial nc plurisubharmonic (nc plush) if it has an nc complex 
hessian that is positive semidefinite when evaluated on all tuples of $n \times n$
matrices for every size $n$; i.e.,
$$
q(X,X^T)[H,H^T] \succeq 0
$$
for all $X,H \in (\mathbb{R}^{n \times n})^g$ for every $n \ge 1$.  

In this paper, we classify all symmetric nc plush polynomials as convex polynomials with 
an nc analytic change of variables; i.e., an nc symmetric polynomial $p$ is nc plush if and only 
if it has the form
\beq \label{eq:ClassifyPlush}
p = \sum f_j^T f_j + \sum k_j k_j^T + F + F^T
\eeq
where the sums are finite and $f_j$, $k_j$, $F$ are all nc analytic.

In this paper, we also present a theory of noncommutative integration 
for nc polynomials and we prove a noncommutative version of the 
Frobenius theorem.

A subsequent paper, \cite{G10}, proves that if the nc complex hessian, $q$, of $p$ 
takes positive semidefinite values on an ``nc open set'' then $q$ takes positive 
semidefinite values on all tuples $X, H$. Thus, $p$ has the form in Equation 
\eqref{eq:ClassifyPlush}. The proof, in \cite{G10}, draws on most of the theorems 
in this paper together with a very different technique involving representations of 
noncommutative quadratic functions.

\end{abstract}

\section{Introduction} \label{sec:Intro}

The main findings of this paper, the description of nc
plush functions and an integration theory for nc functions, 
were indicated in the abstract, so the introduction consists of
precise definitions and precise statements of the plush results. 
The introduction is laid out as follows.

In Section \ref{subsec:BasicDefns}, we review noncommutative (nc) 
polynomials and matrix positivity. In Section \ref{subsec:differentiation}, 
we review nc directional derivatives and provide some examples 
as to how to compute them. We also define the nc complex hessian and 
nc plurisubharmonicity and we compare this to the classical analogue in 
several complex variables. Finally, we close the introduction with the 
statement of the main theorems in Section \ref{subsec:MainResult}. 
These theorems classify all nc plurisubharmonic polynomials and then 
describe the uniqueness of this classification.

\subsection{NC Polynomials, Their Derivatives, and Plurisubharmonicity} \label{subsec:BasicDefns}
Now we give basic definitions. Many of the the definitions we shall 
need sit in the context of an elegant theory of noncommutative analytic 
functions, such as is developed in the articles \cite{KVV} and \cite{Voi1,Voi2}; 
see also \cite{Pop6}. Also, related to our results are those on various classes of noncommutative 
functions on balls as in \cite{AK,BGM}. Transformations on nc variables with 
analytic functions are described in \cite{HKM, Pop6}.

\subsubsection{NC Variables and Monomials} \label{subsubsec:NCVarsMons}
We consider the free semi-group on the $2g$ noncommuting
formal variables $x_1, \ldots, x_g, x_1^T, \ldots, x_g^T$.
The variables $x_j^T$ are the formal transposes of the
variables $x_j$.  The free semi-group in these $2g$
variables generates monomials in all of these variables
$x_1, \ldots, x_g, x_1^T, \ldots, x_g^T$, often called
monomials in $x$, $x^T$. 

If $m$ is a monomial, then
$m^T$ denotes the transpose of the monomial $m$.
For example, given the monomial (in the $x_j$'s)
$x^w=x_{j_1}x_{j_2}\ldots x_{j_n}$,
the involution applied to $x^w$ is
$(x^w)^T=x_{j_n}^T \ldots x_{j_2}^T
x_{j_1}^T$.

\subsubsection{The Ring of NC Polynomials} \label{subsubsec:NCPolys}
Let $\mathbb{R} \langle x_1,\ldots,x_g, x_1^T, \ldots, x_g^T \rangle$ denote
the ring of noncommutative polynomials over $\mathbb{R}$ in the noncommuting
variables $x_1,\ldots,x_g, x_1^T,\ldots, x_g^T$.  We often abbreviate
$$
\mathbb{R} \langle x_1,\ldots,x_g,x_1^T, \ldots, x_g^T \rangle \ \ \
\text{by} \ \ \ \mathbb{R} \langle x,x^T \rangle.
$$
Note that $\mathbb{R}\langle x,x^T \rangle$ maps to itself under the
involution ${}^T$. 

We call a polynomial \textbf{nc analytic} if it contains only the
variables $x_j$ and none of the transposed variables $x_i^T$. Similarly,
we call a polynomial \textbf{nc antianalytic} if it contains only the
variables $x_j^T$ and none of the variables $x_i$.

We call an nc polynomial, $p$, \textbf{symmetric} if $p^T = p$. For
example, $p = x_1 x_1^T + x_2^T x_2$ is symmetric. The
polynomial $\tilde{p} = x_1 x_2 x_4 + x_3 x_1$ is nc analytic but not symmetric.
Finally, the polynomial $\hat{p} = x_2^T x_1^T + 4x_3^T$ is nc antianalytic
but not symmetric.

\subsubsection{Substituting Matrices for NC Variables} \label{subsubsec:SubMats}
If $p$ is an nc polynomial in the variables $x_1,\dots,x_g,x_1^T,\ldots,x_g^T$
and
$$
X=(X_1,X_2,\ldots,X_g) \in (\bbR^{n\times n})^{g},
$$
the evaluation $p(X,X^T)$ is defined by simply replacing $x_j$ by $X_j$ and $x_j^T$ by
$X_j^T$. Note that, for $Z_n=(0_n,0_n,\dots,0_n) \in (\bbR^{n\times n})^{2g}$ where
each $0_n$ is the $n\times n$ zero matrix, $p(0_n)=I_n \otimes p(0_1)$. Because of
this simple relationship, we often simply write $p(0)$ with the size $n$ unspecified.
The involution, ${}^T$, is compatible with matrix transposition, i.e.,
$$
p^T(X,X^T) = p(X,X^T)^T.
$$

\subsubsection{Matrix Positivity} \label{subsubsec:MatrixPositivity}
We say that an nc symmetric polynomial, $p$, in the $2g$ variables 
$x_1,\ldots,x_g,x_1^T,\ldots,x_g^T$, is 
\textbf{matrix positive} if $p(X,X^T)$ is a positive semidefinite matrix when 
evaluated on every $X \in (\mathbb{R}^{n \times n})^{g}$ for every size 
$n \ge 1$; i.e., 
$$
p(X,X^T) \succeq 0
$$
for all $X \in (\mathbb{R}^{n \times n})^{g}$ and all $n \ge 1$.

In \cite{H02}, Helton classified all matrix positive nc symmetric polynomials 
as sums of squares. We recall Theorem 1.1 from \cite{H02}:
\begin{theorem} \label{thm:Helton}
Suppose $p$ is a non-commutative symmetric polynomial. If $p$ is a sum 
of squares, then $p$ is matrix positive. If $p$ is matrix positive, then $p$ 
is a sum of squares.
\end{theorem}

\subsection{NC Differentiation} \label{subsec:differentiation}
First we make some definitions and state some properties about nc differentiation. 
s
In the study of complex variables, we have polynomials in $z$ and $\bar{z}$. 
We can then take derivatives with respect to $z$ and $\bar{z}$; i.e., 
$\frac{\partial p}{\partial z}(z,\bar{z})$ and $\frac{\partial p}{\partial \bar{z}}(z,\bar{z})$. 
We can also make a matrix and fill it with mixed partial derivatives; i.e., the 
$(i,j)$-th entry of the matrix is $\frac{\partial^2 p}{\partial z_i \partial \bar{z}_j}$. 
In classical several complex variables, this matrix of mixed partial derivatives is 
called the complex hessian.

The noncommutative differentiation of polynomials in $x$ and $x^T$ defined in 
this paper is analogous to classical differentiation of polynomials in $z$ and $\bar{z}$ 
from several complex variables.

\subsubsection{Definition of Directional Derivative} \label{subsubsec:dird}
Let $p$ be an nc polynomial in the nc variables $x = (x_1,\ldots, x_g)$ and $x^T = (x_1^T, \ldots, x_g^T)$. 
In order to define a directional derivative, we first replace all $x_i^T$ by $y_i$. Then the \textbf{directional 
derivative of $p$ with respect to $x_j$ in the direction $h_j$} is 
\beq \label{eq:derivx}
p_{x_j}[h_j] := \frac{\partial p}{\partial x_j}(x,x^T)[h_j] = 
\frac{dp}{dt}(x_1, \ldots, x_j+th_j, \ldots, x_g, y_1,\ldots,y_g)|_{t=0} |_{y_i=x_i^T}.
\eeq
The \textbf{directional derivative of $p$ with respect to $x_j^T$ in the direction $k_j$} is 
\beq \label{eq:derivxT}
p_{x_j^T}[k_j] := \frac{\partial p}{\partial x_j^T}(x,x^T)[k_j] = 
\frac{dp}{dt}(x_1,\ldots,x_g, y_1, \ldots, y_j+tk_j, \ldots, y_g)|_{t=0} |_{y_i=x_i^T}.
\eeq
Often, we take $k_j=h_j^T$ in Equation \eqref{eq:derivxT} and we define
\begin{eqnarray*} 
p_x[h] &:=& \frac{\partial p}{\partial x}(x,x^T)[h] = \frac{dp}{dt}(x+th,y)|_{t=0}|_{y=x^T} = 
\sum_{i=1}^{g} {\frac{\partial p}{\partial x_i}(x,x^T)[h_i]} \\
p_{x^T}[h^T] &:=& \frac{\partial p}{\partial x^T}(x,x^T)[h^T] = \frac{dp}{dt}(x,y+tk)|_{t=0}|_{y=x^T, k=h^T} = 
\sum_{i=1}^{g} {\frac{\partial p}{\partial x_i^T}(x,x^T)[h_i^T]}.
\end{eqnarray*}
Then, we (abusively\footnote{For more detail, see \cite{HMV06}. The idea for computing 
$p^{(\ell)}(x)[h]$ is that we first noncommutatively expand $p(x+th)$. Then, 
$p^{(\ell)}(x)[h]$ is the coefficient of $t^\ell$ multiplied by $\ell !$; i.e.,
$
p^{(\ell)}(x)[h] = (\ell !) (\mbox{coefficient of $t^\ell$}).
$ }) 
define the $\ell^{th}$ directional derivative of $p$ in the direction $h$ as 
$$
p^{(\ell)}(x)[h] := \frac{d^{\ell}p}{dt^{\ell}}(x+th,y+tk)|_{t=0}|_{y=x^T,k=h^T}
$$
so the first directional derivative of $p$ in the direction $h$ is 
\begin{eqnarray} \label{eq:dird}
p'(x)]h] &=& \frac{\partial p}{\partial x}(x,x^T)[h] + \frac{\partial p}{\partial x^T}(x,x^T)[h^T] \\
&=& p_x[h] + p_{x^T}[h^T].
\end{eqnarray}
It is important to note that the directional derivative is an nc polynomial that is 
homogeneous degree 1 in $h$, $h^T$.  If $p$ is symmetric, so is $p'$.

\subsubsection{Examples of Differentiation} \label{subsubsec:ex}
Here we provide some examples of how to compute directional derivatives. 

\begin{example} \label{ex:MixedEx}
{\rm Let $p = x_1 x_2^T x_1 + x_1^T x_2 x_1^T$. Then we have 
\begin{eqnarray*}
p_{x_1}[h_1]&=&\frac{\partial p}{\partial x_1}(x,x^T)[h_1] = h_1 x_2^T x_1 + x_1 x_2^T h_1\\
p_{x_2^T}[h_2^T]&=&\frac{\partial p}{\partial x_2^T}(x,x^T)[h_2^T] = x_1 h_2^T x_1\\
p_x[h] &=& \frac{\partial p}{\partial x}(x,x^T)[h] =  h_1 x_2^T x_1 + x_1 x_2^T h_1 + x_1^T h_2 x_1^T 
\end{eqnarray*}
and, 
$$ p'(x)[h] = h_1 x_2^T x_1 + x_1 h_2^T x_1 + x_1 x_2^T h_1 + 
h_1^T x_2 x_1^T + x_1^T h_2 x_1^T +  x_1^T x_2 h_1^T . $$ }
\end{example}

\begin{example} \label{ex:DirDMon}
{\rm Given a general monomial, with $c \in \bbR$,
$$ 
m = c x^{i_1}_{j_1} x^{i_2}_{j_2} \cdots x^{i_n}_{j_n} 
$$
where each $i_k$ is either 1 or $T$, we get that 
$$ 
m' =  c h^{i_1}_{j_1}x^{i_2}_{j_2}\cdots x^{i_n}_{j_n} + 
c x^{i_1}_{j_1}h^{i_2}_{j_2}x^{i_3}_{j_3}\cdots x^{i_n}_{j_n} + \cdots 
+ c x^{i_1}_{j_1}\cdots x^{i_{n-1}}_{j_{n-1}}h^{i_n}_{j_n}. 
$$ }
\end{example}

\subsubsection{Hessian and Complex Hessian} \label{subsubsec:Hessians}
Often, one is most interested in the hessian of a polynomial and its positivity; as 
this determines convexity. However, in this paper, we are most concerned with the 
\textit{complex hessian}, since it turns out to be related to ``nc analytic'' changes 
of variables.

We define the \textbf{nc complex hessian} \index{nc complex hessian}, $q(x,x^T)[h,h^T]$, of
an nc polynomial $p$ as the nc polynomial in the $4g$ variables $x = (x_1,\ldots, x_g)$, $x^T =
(x_1^T, \ldots, x_g^T)$, $h = (h_1, \ldots, h_g)$, and $h^T = (h_1^T, \ldots, h_g^T)$

\beq \label{def:complexhesspoly} \index{$q(x)[h]$}
q(x,x^T)[h,h^T] := \frac{\partial^2 p}{\partial s \partial t}(x+th, y+sk) |_{t,s=0} |_{y=x^T, k=h^T}.
\eeq

The nc complex hessian is an iterated nc directional derivative in the 
sense that we compute it as follows. We first take the nc directional 
derivative of $p$ with respect to $x^T$ in the direction $h^T$ to get 
$p_{x^T}[h^T]$. Then, we take the nc directional derivative of that 
with respect to $x$ in the direction $h$ to get $(p_{x^T}[h^T])_x[h]$. 
We will see later, in Lemma \ref{lem:MixedPartials}, that we can switch 
the order of differentiation to $(p_{x}[h])_{x^T}[h^T]$ and we still 
get the same answer. Sometimes, we note the nc complex hessian as 
$p_{x^T,x}[h^T,h]$. Hence, we have the following equivalent notations 
for the nc complex hessian (and we will use each one when 
context is convenient):
\beq \label{eq:CHnotation}
q(x,x^T)[h,h^T] = p_{x^T,x}[h^T,h] = (p_{x^T}[h^T])_x[h] = (p_{x}[h])_{x^T}[h^T].
\eeq

Note, an extremely important fact about $q(x,x^T)[h,h^T]$, which is restated 
in Theorem \ref{thm:P1andP2} (P1), is that it is quadratic in $h,h^T$ and 
that \textbf{each term contains some $h_j$ and some $h_k^T$}. The nc complex 
hessian is actually a piece of the full \textbf{nc hessian} which is
\begin{eqnarray*}
p'' &=& \frac{\partial^2 p}{\partial t^2}(x+th, y) |_{t=0} |_{y=x^T} +
\frac{\partial^2 p}{\partial t \partial s}(x+th, y+sk) |_{t,s=0} |_{y=x^T, k=h^T} \\
&+& \frac{\partial^2 p}{\partial s \partial t}(x+th, y+sk) |_{t,s=0} |_{y=x^T, k=h^T} +
\frac{\partial^2 p}{\partial s^2}(x, y+sk) |_{s=0} |_{y=x^T, k=h^T} \\
&=& 2q(x,x^T)[h,h^T] + \frac{\partial^2 p}{\partial t^2}(x+th, y) |_{t=0} |_{y=x^T} +
\frac{\partial^2 p}{\partial s^2}(x, y+sk) |_{s=0} |_{y=x^T, k=h^T}.
\end{eqnarray*}

We call a symmetric nc polynomial, $p$, \textbf{nc plurisubharmonic}
(or \textbf{nc plush}) if the nc complex hessian, $q$, of $p$ is matrix 
positive. In other words, we require that $q$ be positive semidefinite 
when evaluated on all tuples of $n \times n$ matrices for every
size $n$; i.e.,
$$
q(X,X^T)[H,H^T] \succeq 0
$$
for all $X,H \in \gtupn$ for every $n \ge 1$. 

\subsubsection{Examples of Complex Hessians} \label{subsubsec:ExCH}
Here we provide some examples of how to compute nc 
complex hessians.

\begin{example} \label{ex:MixedExCH}
{\rm Let $p = x_1 x_2^T x_1 + x_1^T x_2 x_1^T$ as in Example \ref{ex:MixedEx}. 
Then, we have
$$
q = h_1 h_2^T x_1 + x_1 h_2^T h_1 + h_1^T h_2 x_1^T + x_1^T h_2 h_1^T.
$$ }
\end{example}

\begin{example} \label{ex:ExCH}
{\rm Let $p = x^T x^T x x$. Then, we have 
\begin{eqnarray*}
q(x,x^T)[h,h^T] &=& h^T x^T h x + h^T x^T x h + x^T h^T h x + x^T h^T x h \\
 &=& (h^T x^T + x^T h^T)(hx + xh) \\
 &=& (hx + xh)^T (hx + xh).
\end{eqnarray*}
We can see that, for any $X,H \in \mathbb{R}^{n\times n}$ for any size 
$n \ge 1$, we have that
$$
q(X,X^T)[H,H^T] = (HX + XH)^T (HX + XH) \succeq 0.
$$
Hence, this nc polynomial, $p = x^T x^T x x$, is nc plush.}
\end{example}

\begin{example} \label{ex:AntiAnaCH}
{\rm The nc complex hessian of any nc analytic polynomial is 0. 
The nc complex hessian of any nc antianalytic polynomial is 0. 
Hence, both nc analytic and nc antianalytic polynomials are 
nc plush. }
\end{example}

\subsection{Main Result} \label{subsec:MainResult}
In this paper we classify all symmetric nc plush polynomials in $g$ free variables. 

\begin{theorem} \label{thm:MainThm}
An nc symmetric polynomial $p$ in free variables is nc plurisubharmonic 
if and only if $p$ can be written in the form
\beq \label{eq:MainEqn}
p = \sum f_j^T f_j + \sum k_j k_j^T + F + F^T
\eeq
where the sums are finite and each $f_j$, $k_j$, $F$ is nc analytic.
\end{theorem}

\proof The proof requires the rest of this paper and culminates in 
Section \ref{sec:ProofMainResult}. \qed 

The subsequent paper (\cite{G10}) strengthens the result of Theorem 
\ref{thm:MainThm} by weakening the hypothesis while keeping the 
same conclusion. Specifically, we assume that the nc polynomial 
is nc plush on ``an nc open set'' and conclude that it is nc plush 
everywhere and hence has the form in Equation \eqref{eq:MainEqn}. 
The proof, in \cite{G10}, draws on most of the theorems in this paper 
together with a very different technique involving representations of 
noncommutative quadratic functions.

The representation in Equation \eqref{eq:MainEqn} is unique up to 
the natural transformations.

\begin{theorem} \label{thm:MainThmUnique}
Let $p$ be an nc symmetric polynomial in free variables that is nc plurisubharmonic 
and let
\begin{eqnarray*}
\tN &:=& \min \{N : p = \sum_{j=1}^{N}{f_j^T f_j} + \sum_{j=1}^{M}{k_j k_j^T} + F + F^T\} \\
\tM &:=& \min \{M: p = \sum_{j=1}^{N}{f_j^T f_j} + \sum_{j=1}^{M}{k_j k_j^T} + F + F^T\}.
\end{eqnarray*}
Then, we can represent $p$ as 
$$
p = \sum_{j=1}^{\tN}{\tf_j^T \tf_j} + \sum_{j=1}^{\tM}{\tk_j \tk_j^T} + \tilde{F} +\tilde{F}^T
$$
and if $N$ and $M$ are integers such that $N \ge \tN$, $M \ge \tM$ and 
$$
p = \sum_{j=1}^{N}{f_j^T f_j} + \sum_{j=1}^{M}{k_j k_j^T} + F + F^T \ ,
$$
then there exist isometries $U_1 : \mathbb{R}^{\tN} \longrightarrow \mathbb{R}^{N}$ and 
$U_2 : \mathbb{R}^{\tM} \longrightarrow \mathbb{R}^{M}$ such that 
$$
\left( \begin{array}{c}
f_1 \\
\vdots \\
f_N
\end{array} \right)
=
U_1 
\left( \begin{array}{c}
\tf_1 \\
\vdots \\
\tf_{\tN}
\end{array} \right)
+ 
\vec{c}_1
\ \ \ \mbox{and} \ \ \ 
\left( \begin{array}{c}
k_1^T \\
\vdots \\
k_M^T
\end{array} \right)
=
U_2
\left( \begin{array}{c}
\tk_1^T \\
\vdots \\
\tk_{\tM}^T
\end{array} \right)
+ 
\vec{c}_2
$$
where $\vec{c}_1 \in \mathbb{R}^N$ and $\vec{c}_2 \in \mathbb{R}^M$.
\end{theorem}

\proof Theorem \ref{thm:MainThm} gives the desired form of $p$ and 
nc integration will give the uniqueness.  We provide the details of the 
proof in Section \ref{sec:ProofMainResult}. \qed 

A byproduct of the proof of Theorem \ref{thm:MainThm} is noncommutative 
integration theory of nc polynomials. This includes a Frobenius theorem for 
nc polynomials and is discussed further in Section \ref{sec:integration}.

\subsection{Guide to the Paper} \label{subsec:Guide}
In Section \ref{sec:integration}, we provide a theory of noncommutative 
integration for nc polynomials and in Section \ref{subsec:NCIntegration}, 
we state and prove a noncommutative version of the Frobenius theorem. 
In Section \ref{sec:ComplexHessSOS}, we prove that the nc complex hessian 
for an nc plush polynomial is the sum of hereditary and antihereditary squares. 
Finally, in Section \ref{sec:ProofMainResult}, we prove the main results. We 
apply nc integration theory to the sum of squares representation of 
the nc complex hessian found in Section \ref{sec:ComplexHessSOS}. 
We also settle the issue of uniqueness of this sum of squares representation.

The authors would like to thank Mark Stankus for fruitful discussions on 
noncommutative integration.

\section{NC Integration} \label{sec:integration}

In this section, we introduce a natural notion of noncommutative (nc) 
integration and then give some basic properties. We say that an nc polynomial 
$p$ in $x=(x_1,\ldots, x_g)$ and $h_j$ is \textbf{integrable in $x_j$} 
if there exists an nc polynomial $f(x)$ such that $f_{x_j}[h_j] = p$. 
We say that an nc polynomial $p$ in $x=(x_1,\ldots, x_g)$ and 
$h=(h_1,\ldots, h_g)$ is \textbf{integrable} if there exists an nc 
polynomial $f(x)$ such that $f'(x)[h] = p$.

\subsection{Notation} \label{subsec:Notation}

Let $m$ be a monomial containing only the variables $x_1$, $x_2$, $\ldots$, 
$x_g$. When we write $\monsub{m}{x_i}{h_i}$, we mean the set of monomials 
that are degree one in $h_i$ where one $x_i$ in $m$ has been replaced by $h_i$. For 
example, if $m = x_1 x_2 x_1 x_2$, then 
$$
\monsub{m}{x_1}{h_1} = \{ h_1 x_2 x_1 x_2, \ x_1 x_2 h_1 x_2 \}
$$
and 
$$
\monsub{m}{x_2}{h_2} = \{ x_1 h_2 x_1 x_2, \ x_1 x_2 x_1 h_2 \}.
$$ 
We also define a double substitution as follows. When we write 
$$
\monsubtwo{m}{x_i}{h_i}{x_j}{h_j} := \monsub{ \left( \monsub{m}{x_i}{h_i} \right) }{x_j}{h_j},
$$ 
we mean 
the set of monomials that are degree one in $h_i$ and degree one in 
$h_j$ where one $x_i$ in $m$ has been replaced by $h_i$ and 
one $x_j$ in $m$ has been replaced by $h_j$. Note that we have 
\beq \label{eq:TwoMonSubEq}
\monsubtwo{m}{x_i}{h_i}{x_j}{h_j} = \monsubtwo{m}{x_j}{h_j}{x_i}{h_i}.
\eeq

Using $m = x_1 x_2 x_1 x_2$, we have that 
\begin{eqnarray*}
\monsubtwo{m}{x_1}{h_1}{x_2}{h_2} &=& \monsubtwo{m}{x_2}{h_2}{x_1}{h_1} \\
&=& \{h_1 h_2 x_1 x_2, \ h_1 x_2 x_1 h_2,  \ x_1 h_2 h_1 x_2, \ 
x_1 x_2 h_1 h_2 \}.
\end{eqnarray*}

Sometimes we will start with a monomial $m$ that is degree 1 in $h_i$
and we wish to replace this $h_i$ by $x_i$. When we write $\monsub{m}{h_i}{x_i}$, 
the set we get contains just one monomial so we abuse notation and 
use $\monsub{m}{h_i}{x_i}$ to represent the actual monomial in this 
set.

\subsection{Differentially Wed Monomials} \label{subsec:Wed}

For $\gamma$ either 1 or $T$, two monomials $m$ and $\tm$ are called
{\bf 1-differentially wed with respect to $x^\gamma_j$} if both 
$m$ and $\tm$ contain the same $h^\gamma_j$ to exactly first order
and if $m$ has an $x^\gamma_j$ where $\tm$ has an $h^\gamma_j$
and if $\tm$ has an $x^\gamma_j$ where $m$ has an $h^\gamma_j$.
Thus interchanging $h^\gamma_j$ and this $x^\gamma_j$ in $m$
produces $\tm$; i.e.,
$$
\monsub{m}{h_j^\gamma}{x_j^\gamma} = \monsub{\tm}{h_j^\gamma}{x_j^\gamma}
$$

More generally, for $\alpha,\beta$ either 1 or $T$, two monomials $m$ 
and $\tm$ are called \textbf{1-differentially wed} if both are degree 
one in $h$ or $h^T$ and if 
$$
\monsub{m}{h_i^\alpha}{x_i^\alpha} = \monsub{\tm}{h_j^\beta}{x_j^\beta}
$$

From these definitions, if $m$ and $\tm$ are 1-differentially 
wed with respect to a particular variable then $m$ and $\tm$ 
are 1-differentially wed but not the other way around (which 
we demonstrate below). 

\begin{example} \label{ex:DiffWed1}
{\rm The monomials $m = h_1 x_2^T x_1$ and $\tm = x_1 x_2^T h_1$ are 
1-differentially wed with respect to $x_1$. }
\end{example}

\begin{example} \label{ex:DiffWed2}
{\rm The monomials $m = h_1 x_2^T x_1$ and $\tm = x_1 h_2^T x_1$ are 
1-differentially wed (but not with respect to a particular variable).}
\end{example}

\begin{example} \label{ex:DiffWed3}
{\rm The monomials $m = x_2 h_2 x_2$ and $\tm = x_1 x_2 h_2$ are not 
1-differentially wed with respect to any variable.}
\end{example}
 
\begin{theorem} \label{thm:intGvars}
   A polynomial $p$ in $x=(x_1,\ldots,x_g),h=(h_1,\ldots,h_g)$ 
   is integrable if and 
   only if each monomial in $p$ has degree one in $h$ (i.e., 
   contains some $h_j$) and whenever a monomial $m$ 
   occurs in $p$, each monomial which is 1-differentially wed 
   to $m$ also occurs in $p$ and has the same coefficient.
\end{theorem}

\proof First suppose the polynomial $p$ in $x, h$ can be integrated 
in $x$. Then there exists a polynomial, $f(x)$, such that $f'(x)[h] = p$. 
Write $f$ as 
$$
f = \sum_{i=1}^{N}{m_i}
$$
where each $m_i$ is a monomial in $x$. Then, by applying Example 
\ref{ex:DirDMon} to each monomial $m_i$, if a monomial $\tm$ occurs 
in $p = f'$, then every 1-differentially wed monomial to $\tm$ also 
occurs in $p=f'$ with the same coefficient.

Now suppose each monomial in $p$ has degree 1 in $h$ (i.e., contains 
some $h_j$) and if $m$ is a monomial in $p$, then each monomial which 
is 1-differentially wed also occurs in $p$ with the same coefficient.  
We will show that $p$ can be integrated in $x$

Write 
$$ p = \sum_{i=1}^{N} {m_i} $$
where $m_i$ is a monomial in $x$ and degree 1 in $h$. 
Now we will change the order of summation of these monomials 
so that we group together all monomials that are 1-differentially wed. 
We do this in the following way. 

Let $w_1$ be the polynomial that contains $m_1$ and all 
1-differentially wed monomials to $m_1$. 

Let $1 \le \alpha_2$ be the smallest integer such that $m_{\alpha_2}$ 
is not a term in $w_1$. Then let $w_2$ be the polynomial 
that contains $m_{\alpha_2}$ and all 1-differentially wed monomials to 
$m_{\alpha_2}$.

Let $1 \le \alpha_3$ be the smallest integer such that $m_{\alpha_3}$ 
is not a term in $w_1$ and not a term in $w_2$. Then let 
$w_3$ be the polynomial that contains $m_{\alpha_3}$ and 
all 1-differentially wed monomials to $m_{\alpha_3}$.

We continue this process until it stops (it stops since $p$ is a finite sum of monomials). 
Then we have written $p$ as 
$$ p = \sum_{i=1}^{\ell} w_i. $$

It is important to note that with this construction, each $w_i$ is a 
homogeneous polynomial of some fixed degree where each monomial 
in $w_i$ is degree 1 in $h$. 

Now define $\alpha_1 = 1$ and
$$
f_i (x) := \monsub{m_{\alpha_i}}{h}{x} , \ \ \  1 \le i \le \ell.
$$

Then we have, by properties of differentiation and construction of $w_i$, 
that $f_i' = w_i$.  Finally, define 
$$ f(x) := \sum_{i=1}^{\ell} f_i (x) $$
and notice that $f' = p$.   \qed  \\

\begin{cor}
\label{cor:int1var}
   A polynomial $p$ in $x,h_j$ is integrable in $x_j$ 
   if and only if 
   each monomial in $p$ has degree one in $h_j$ and
   whenever a monomial $m$ occurs in $p$, each monomial
   which is 1-differentially wed with respect to $x_j$ 
   also occurs in $p$ and has the same coefficient.
\end{cor}

\subsection{Uniqueness of Noncommutative Integration} \label{subsec:UniqueInteg}
In this subsection, we explore the uniqueness of noncommutative integration. 
In classical calculus, integrating produces constants of integration. Here, we 
provide the noncommutative analogue.

\begin{prop} \label{prop:MonDerivUnique}
Suppose $m$ and $\tm$ are distinct monomials in the variables 
$x=(x_1,\ldots, x_g)$. Then, we have that
\begin{enumerate}
\item $(m)_{x_i}[h_i]$ and $(\tm)_{x_i}[h_i]$ have no terms in common and 
hence 
$$
(m)_{x_i}[h_i] \ne (\tm)_{x_i}[h_i]
$$
provided $x_i$ is contained in either $m$ or $\tm$; and
\item we have that 
$$
m' = (m)_{x}[h] \ne (\tm)_{x}[h] = \tm'.
$$
\end{enumerate}
If $m$ and $\tm$ are distinct monomials in the variables $x=(x_1, \ldots, x_g)$ 
and $y=(y_1, \ldots, y_s)$, then 
$$
m_x[h] \ne \tm_x[h].
$$
\end{prop}

\proof If $m$ and $\tm$ have different degree, then so do 
their derivatives and we are done. Suppose $m$ and $\tm$ have 
the same degree. If $m$ contains $x_i$ and $\tm$ does not, 
then $(\tm)_{x_i}[h_i] = 0$ while $(m)_{x_i}[h_i]$ is a nonzero nc 
polynomial.

Suppose both $m$ and $\tm$ contain $x_i$ and are the same degree. 
Then, write 
$$
m = x_{j_1} x_{j_2} \cdots x_{j_s} \ \ \ \ \ \mbox{and} \ \ \ \ \  \tm = x_{k_1} x_{k_2} \cdots x_{k_s} 
$$
where the tuple of integers $(j_1, j_2, \ldots, j_s)$ is not 
the same as the tuple of integers $(k_1, k_2, \ldots, k_s)$. 
This forces $(m)_{x_i}[h_i]$ and $(\tm)_{x_i}[h_i]$ to have 
no terms in common. This completes the proof of \textit{(1)}.

To prove \textit{(2)}, note that 
$$
m' = (m)_x[h] = \sum_{i=1}^{g}{(m)_{x_i}[h_i]} \ \ \ \ \ \mbox{and} \ \ \ \ \
\tm' = (\tm)_x[h] = \sum_{i=1}^{g}{(\tm)_{x_i}[h_i]}
$$
and if $m' = \tm'$, then we must have $(m)_{x_i}[h_i] = \tm_{x_i}[h_i]$ 
for each $i$. However, \textit{(1)} implies that this is impossible. 

If $m$ and $\tm$ are distinct monomials in the variables $x=(x_1,\ldots,x_g)$ 
and $y=(y_1,\ldots,y_g)$, the proof follows exactly the way the proof of 
\textit{(2)} does. \qed 

\begin{lemma}  \label{lem:derivPolyOne}
Suppose $p$ is an nc polynomial in the variables $x=(x_1,\ldots, x_g)$ such that  
$p_{x_i}[h_i]=0$. Then, $p(x_1,\ldots,x_g) = f(x_1,\ldots, x_{i-1}, x_{i+1},\ldots, x_g)$ 
is an nc polynomial in the variables $x_1,\ldots, x_{i-1},x_{i+1},\ldots, x_g$.
\end{lemma}

\proof First, if $m$ is a monomial in the variables $x=(x_1,\ldots,x_g)$ 
that contains $x_i$, then $m_{x_i}[h_i]$ is 
a sum of terms where each instance of $x_i$ is replaced by $h_i$ (see 
Example \ref{ex:MixedEx}). Note that each term in $m_{x_i}[h_i]$ has a 
different number of variables to the left of $h_i$; hence, the terms can 
not cancel. Thus, $m_{x_i}[h_i] \ne 0$.

Now suppose $p$ is an nc polynomial in the variables $x=(x_1,\ldots,x_g)$. 
We write the nc polynomial $p$ as
\beq \label{eq:ind}
p = \sum_{j=1}^{N} {\alpha_j \, m_j}
\eeq
where the $\alpha_j$ are nonzero real constants and the $m_j$ are distinct monomials. 
Then, we have that 
\beq \label{eq:cancel}
0 = p_{x_i}[h_i] = \sum_{j=1}^{N}{\alpha_j \, (m_j)_{x_i}[h_i]}.
\eeq
Since the $m_j$ are distinct monomials, Proposition \ref{prop:MonDerivUnique} 
implies that no cancellation can occur in Equation \eqref{eq:cancel}. 
This implies that 
$$
(m_j)_{x_i}[h_i] = 0
$$
for all $j = 1, \ldots, N$. Then, by the first paragraph in this proof, we get that each 
$m_j$ is a monomial in the variables $x_1, \ldots, x_{i-1},x_{i+1},\ldots,x_g$. 
This implies that $p$, as in Equation \eqref{eq:ind}, is a polynomial in 
the variables $x_1, \ldots, x_{i-1}, x_{i+1}, \ldots, x_g$. \qed 

\begin{prop} \label{prop:ZeroDeriv}
Suppose $p$ is an nc polynomial in the $g+s$ variables 
$x=(x_1, \ldots, x_g)$ and $y=(y_1, \ldots, y_s)$. If 
$p_x[h] = 0$, then 
$$
p(x,y) = p(x_1,\ldots,x_g,y_1,\ldots,y_s) = f(y_1,\ldots,y_s)
$$
is an nc polynomial in the variables $y=(y_1, \ldots, y_s)$.
\end{prop}

\proof If $p$ is an nc polynomial in the $g+s$ variables 
$x = (x_1, \ldots, x_g)$ and $y = (y_1, \ldots, y_s)$, then 
$$
p_x[h] = \sum_{i=1}^{g}{p_{x_i}[h_i]}.
$$
It is important to note that $p_{x_i}[h_i]$ is an nc polynomial 
in $x$, $y$, and $h_i$. Since 
$$
p_x[h] = \sum_{i=1}^{g}{p_{x_i}[h_i]} = 0
$$
and since each $p_{x_i}[h_i]$ is an nc polynomial that is 
linear in $h_i$, it must follow that
$$
p_{x_i}[h_i] = 0 \ \ \ \ \ \forall \ i=1,\ldots, g.
$$
Then, Lemma \ref{lem:derivPolyOne} implies that 
$p$ is an nc polynomial in the variables 
$$
x_1, \ldots, x_{i-1}, x_{i+1}, \ldots, x_g, y_1, \ldots, y_s
$$
for all $1 \le i \le g$. This can only happen if $p$ is 
an nc polynomial in the variables $y = (y_1, \ldots, y_s)$. \qed 

\begin{cor} \label{cor:Constant}
Suppose $p$ is an nc polynomial in $x=(x_1,\ldots,x_g)$. Then, 
we have that
\begin{enumerate}
\item if 
$$
p'(x)[h] = p_x[h] = \sum_{i=1}^{g}{p_{x_i}[h_i] }= 0,
$$
then $p$ is constant, and
\item if $\tp$ is another nc polynomial in the variables 
$x=(x_1,\ldots, x_g)$ such that $p' = \tp'$ then $p = \tp + \alpha$ 
where $\alpha$ is a real constant.
\end{enumerate}
\end{cor}

\proof Property \textit{(1)} directly follows from Proposition \ref{prop:ZeroDeriv}.

If $p' = \tp'$, then we get that
$$
0 = p' - \tp' = (p - \tp)'
$$
which, by property \textit{(1)}, implies that $p-\tp$ is constant. \qed 

\subsubsection{Noncommutative Complex Differentiation} \label{subsubsec:NCCompDiff}
Here, we specialize from the variables $x=(x_1,\ldots, x_g)$ and $y=(y_1,\ldots,y_s)$ 
to the variables $x=(x_1,\ldots, x_g)$ and $x^T=(x_1^T,\ldots,x_g^T)$. The 
first corollary below, Corollary \ref{cor:AA}, follows directly from Proposition 
\ref{prop:ZeroDeriv} above.

\begin{cor} \label{cor:AA}
Suppose $p$ is an nc polynomial in the variables $x=(x_1,\ldots,x_g)$ and 
$x^T = (x_1^T, \ldots, x_g^T)$. Then, we have that
\begin{enumerate}
\item if $p_x[h] = 0$, then $p$ is an nc antianalytic polynomial, and
\item if $p_{x^T}[h^T] = 0$, then $p$ is an nc analytic polynomial.
\end{enumerate}
\end{cor}

\begin{lemma}\label{lem:leviUnique}
Let $p$ be an nc polynomial in the nc variables $x=(x_1,\ldots, x_g)$, 
$x^T = (x_1^T, \ldots, x_g^T)$ and let $q$ be the nc complex hessian 
of $p$.  Then $q = 0$ if and only if $p = F + G^T$ where $F$ and $G$ 
are nc analytic polynomials.

If, in addition, $p$ is symmetric, then $q = 0$ if and only if $p = F + F^T$ 
where $F$ is an nc analytic polynomial.
\end{lemma}

\proof Lemma \ref{lem:MixedPartials}, in Section \ref{subsec:NCIntegration} 
below, allows us to switch the order of differentiation to get
$$
q = p_{x^T,x}[h^T,h] = (p_{x^T}[h^T])_x[h] = (p_{x}[h])_{x^T}[h^T].
$$
Then, we get that the nc complex hessian of $p = F + G^T$ is
$$
q = (F_{x^T}[h^T])_{x}[h] + (G^T_{x}[h])_{x^T}[h^T] = 0.
$$

Now suppose $p$ contains a term with both $x$ and $x^T$. 
Write $p$ as 
$$
p = \sum_{j=1}^{N} {\alpha_j m_j}
$$
where $\alpha_j$ are nonzero real constants and $m_j$ are 
distinct monomials in $x$ and/or $x^T$. Then, the nc complex 
hessian of $p$ is 
$$
q = \sum_{j=1}^{N} {\alpha_j (m_j)_{x^T,x}[h^T,h] }.
$$
Since the $m_j$ are distinct, Proposition \ref{prop:MonDerivUnique} 
implies that the nc polynomials $(m_i)_{x^T}[h^T]$ and 
$(m_j)_{x^T}[h^T]$ have no terms in common for all $i \ne j$. 
Then, we apply Proposition \ref{prop:MonDerivUnique} again 
to get that the nc polynomials $(m_i)_{x^T,x}[h^T,h]$ and 
$(m_j)_{x^T,x}[h^T,h]$ have no terms in common for all $i \ne j$. 
This implies that no cancellation occurs in $q$ so that $q \ne 0$.  \qed 

\subsection{NC ``Gradient'' of a Potential} \label{subsec:NCIntegration}

In this subsection, we give a noncommutative Frobenius Theorem and 
present some equivalent tests to determine if a list of nc polynomials is 
simultaneously integrable. 

\begin{lemma} \label{lem:MixedPartials}
Suppose $p(x_1,\ldots, x_g)$ is an nc polynomial in $x_1,\ldots, x_g$. 
Then 
$$
(p_{x_i}[h_i])_{x_j}[h_j] = (p_{x_j}[h_j])_{x_i}[h_i].
$$
\end{lemma}

\proof Write 
$$
p = \sum_{\alpha=1}^{t} {m_\alpha}
$$
where each $m_\alpha$ is a monomial in $x_1, \ldots, x_g$. 
Then we have 
$$
(p_{x_i}[h_i])_{x_j}[h_j] = \sum_{\alpha=1}^{t} {((m_\alpha)_{x_i}[h_i])_{x_j}[h_j]}
$$
and
$$
(p_{x_j}[h_j])_{x_i}[h_i] = \sum_{\alpha=1}^{t} {((m_\alpha)_{x_j}[h_j])_{x_i}[h_i]}.
$$
Note that $(m_\alpha)_{x_i}[h_i]$ is 
the sum of all monomials in the set $\monsub{m_\alpha}{x_i}{h_i}$ 
and $((m_\alpha)_{x_i}[h_i])_{x_j}[h_j]$ is the sum of all monomials 
in the set $\monsubtwo{m_\alpha}{x_i}{h_i}{x_j}{h_j}$. Equation 
\eqref{eq:TwoMonSubEq} implies that 
$$
\monsubtwo{m_\alpha}{x_i}{h_i}{x_j}{h_j} = \monsubtwo{m_\alpha}{x_j}{h_j}{x_i}{h_i}
$$
which implies that 
$$
((m_\alpha)_{x_i}[h_i])_{x_j}[h_j] = ((m_\alpha)_{x_j}[h_j])_{x_i}[h_i].
$$
Hence, $(p_{x_i}[h_i])_{x_j}[h_j] = (p_{x_j}[h_j])_{x_i}[h_i]$. \qed 

The following theorem is the noncommutative analogue of the Frobenius Theorem 
in that the classical specialization of $(a) \Leftrightarrow (b)$ to $x \in \mathbb{R}^g$ 
in the theorem says that 
$$
\left(
\begin{array}{cccc}
f_1, & f_2, & \ldots, & f_g
\end{array}
\right)
$$
is the gradient of a function if and only if 
$$
\frac{\partial f_i}{\partial x_j} = \frac{\partial f_j}{\partial x_i}.
$$

\begin{theorem} \label{thm:NCFrob}
Suppose $\delta$ is an nc polynomial such that 
$$
\delta(x_1,\ldots,x_g, h_1,\ldots, h_g) = \sum_{i=1}^{g} {f_i(x_1,\ldots,x_g,h_i)}
$$
where each $f_i(x_1,\ldots, x_g,h_i)$ is homogeneous of degree 1 in $h_i$. 
Then the following are equivalent:
\begin{itemize}
\item[(a)] $\delta$ is integrable.
\item[(b)] Each $f_i(x_1,\ldots,x_g,h_i)$ is integrable in $x_i$ and 
$ (f_i)_{x_j}[h_j] = (f_j)_{x_i}[h_i]$ for any $i,j$.
\item[(c)] For each monomial, $m$, in $\delta$, every 1-differentially wed 
monomial to $m$ also occurs in $\delta$ (with the same coefficient).
\end{itemize}
\end{theorem}

\proof Theorem \ref{thm:intGvars} gives the equivalence of (a) and (c).

Now we show (a) and (b) are equivalent. First suppose (a) holds. Then 
there exists an nc polynomial $\cP(x_1,\ldots,x_g)$ such that 
$$
\cP' = \delta \Longrightarrow \sum_{i=1}^{g}{\cP_{x_i}[h_i]} = \sum_{i=1}^{g}{f_i(x_1,\ldots,x_g,h_i)}.
$$
This forces $\cP_{x_i}[h_i] = f_i(x_1,\ldots,x_g,h_i)$ and then 
Lemma \ref{lem:MixedPartials} gives that 
$$
(f_i)_{x_j}[h_j] = \left( \cP_{x_i}[h_i] \right)_{x_j} [h_j] 
=  \left( \cP_{x_j}[h_j] \right)_{x_i} [h_i] = (f_j)_{x_i}[h_i]. 
$$

Now suppose $\delta$ is not integrable. Then there exists some monomial $m$ in $\delta$ 
such that not all 1-differentially wed monomials to $m$ occur in $\delta$. Without 
loss of generality, suppose $m$ is a term in $f_1(x,h_1)$. Recall, this implies 
$m$ is degree one in $h_1$.

If $\delta$ does not contain a monomial that is 
1-differentially wed to $m$ with respect to $x_1$, then $f_1(x,h_1)$ is 
not integrable in $x_1$.

Suppose $m$ contains the variable $x_k$ and that $\delta$ (more specifically, $f_k(x,h_k)$) 
does not contain the monomial $\tm$, a specific monomial in the set 
$\monsub{ \left( \monsub{m}{h_1}{x_1} \right) }{x_k}{h_k} 
= \monsub{ \left( \monsub{m}{x_k}{h_k} \right) }{h_1}{x_1}$. Note that $\tm$ is 
1-differentially wed to $m$ and this implies that the sets $\monsub{m}{x_k}{h_k}$ 
and $\monsub{\tm}{x_1}{h_1}$ are equal.
If $(f_1)_{x_k}[h_k] = (f_k)_{x_1}[h_1]$, then the monomial $\hat{\tm}$, a specific 
monomial in the set $\monsub{m}{x_k}{h_k}=\monsub{\tm}{x_1}{h_1}$, is a term in 
$(f_1)_{x_k}[h_k] = (f_k)_{x_1}[h_1]$. This implies that $\tm = \monsub{\hat{\tm}}{h_1}{x_1}$ 
is a term in $f_k(x,h_k)$ which is contained in $\delta$.

Thus, we have shown that if $\delta$ is not integrable then either some $f_i(x,h_i)$ is not 
integrable with respect to $x_i$ or $(f_i)_{x_j}[h_j] \ne (f_j)_{x_i}[h_i]$ for some 
$i \ne j$. \qed 

\subsection{Levi-differentially Wed Monomials} \label{subsec:LeviDiffWed}

Now we turn to properties of the nc complex hessian $q$, as $q$ 
is just a second nc directional derivative. 

Two monomials $m$ and $\tm$ are called \textbf{Levi-differentially wed} 
if $m,\tm$ are both degree 2 in $h,h^T$, $m$ contains some $h_i$, $h_j^T$, 
$\tm$ contains some $h_k$, $h_s^T$ and
$$
\monsubtwo{m}{h_i}{x_i}{h_j^T}{x_j^T} = \monsubtwo{\tm}{h_k}{x_k}{h_s^T}{x_s^T}
$$

Indeed, Levi-differentially wed is an equivalence relation on the monomials 
in the nc complex hessian, $q$, and the coefficients of all Levi-differentially 
wed monomials in $q$ are the same.

\begin{example} \label{ex:LeviWed1}
{\rm The monomials $h^T h x^T x$, $h^T x x^T h$, $x^T h h^T x$, and $x^T x h^T h$ 
are all Levi-differentially wed to each other. }
\end{example}

\begin{example} \label{ex:LeviWed2}
{\rm None of the monomials $h^T h x^T x$, $h^T x h^T x$, $x^T h x^T h$ are 
Levi-differentially wed to each other.}
\end{example}

\begin{theorem} \label{thm:P1andP2}
An nc polynomial $q$ in $x,x^T,h,h^T$ is an nc complex hessian if and only 
if the following two conditions hold:

\begin{itemize}
\item[(P1)] Each monomial in $q$ contains exactly one $h_j$ 
and one $h_k^T$ for some $j,k$. 

\item[(P2)] If a certain monomial $m$ is contained in $q$, any 
monomial $\tm$ that is Levi-differentially wed to $m$ is also 
contained in $q$.
\end{itemize}
\end{theorem}

\proof First suppose $q$ is an nc complex hessian. Equation \eqref{def:complexhesspoly} 
shows that $q$ is an nc directional derivative of an nc directional derivative. Then 
properties of nc directional derivatives imply that (P1) and (P2) hold.

Now suppose (P1) and (P2) hold. Write $q$ as 
$$
q = \sum_{i=1}^{N} {m_i}
$$
where each $m_i$ is a monomial that is quadratic in $h,h^T$ and 
contains some $h_j$ and $h^T_k$.  Now we will change the order of 
summation so that we group together all monomials that are 
Levi-differentially wed to each other. We do this in the following way.

Let $w_1$ be the nc polynomial that contains $m_1$ and all 
Levi-differentially wed monomials to $m_1$. 

Let $1 \le \alpha_2$ be the smallest integer such that $m_{\alpha_2}$ 
is not a term in $w_1$. Then let $w_2$ be the nc 
polynomial that contains $m_{\alpha_2}$ and all Levi-differentially 
wed monomials to $m_{\alpha_2}$.

Let $1 \le \alpha_3$ be the smallest integer such that $m_{\alpha_3}$ 
is not a term in $w_1$ and not a term in $w_2$. Then 
let $w_3$ be the nc polynomial that contains $m_{\alpha_3}$ 
and all Levi-differentially wed monomials to $m_{\alpha_3}$.

We continue this process until it stops (it stops since $q$ is a finite 
sum of monomials). Then we have written $q$ as 
$$
q = \sum_{i=1}^{\ell} {w_i}.
$$
It is important to note that with this construction, each $w_i$ 
is a homogeneous polynomial of some fixed degree where each monomial 
in $w_i$ contains some $h_j$ and some $h_k^T$.

Now define $\alpha_1 = 1$ and
$$
f_i(x,x^T) := \monsubtwo{m_{\alpha_i}}{h}{x}{h^T}{x^T}, \ \ \ 1 \le i \le \ell.
$$
Then we have, by properties of differentiation and construction of 
$w_i$, that the nc complex hessian of each $f_i$ is just 
$w_i$. Finally, define
$$
f(x,x^T) := \sum_{i=1}^{\ell} {f_i(x,x^T)}
$$
and notice that the nc complex hessian of $f$ is $q$. \qed 

\begin{lemma}\label{lem:monwed}
Let $m,m',n,n'$ all be nc analytic monomials with degree 1 in $h$ 
(or all nc antianalytic monomials with degree 1 in $h^T$). Then 
$m$, $m'$ are 1-differentially wed and $n$, $n'$ are 1-differentially wed 
if and only if $n^T m$ and $n'^T m'$ are Levi differentially wed.
\end{lemma}

\proof Without loss of generality, suppose $m,m',n,n'$ are all nc analytic 
monomials with degree 1 in $h$.

$m,m'$ are 1-differentially wed and $n,n'$ are 1-differentially wed 
if and only if 
\begin{eqnarray*}
\monsub{m}{h_i}{x_i} &=& \monsub{m'}{h_j}{x_j} \\
\monsub{n}{h_k}{x_k} &=& \monsub{n'}{h_s}{x_s}.
\end{eqnarray*}
This happens if and only if 
\begin{eqnarray*}
\monsubtwo{(n^T m)}{h^T_k}{x^T_k}{h_i}{x_i} &=& 
\left( \monsub{n}{h_k}{x_k} \right)^T \left( \monsub{m}{h_i}{x_i} \right) \\
&=&  \left( \monsub{n'}{h_s}{x_s} \right)^T \left( \monsub{m'}{h_j}{x_j} \right) = 
\monsubtwo{(n'^T m')}{h^T_s}{x^T_s}{h_j}{x_j}.
\end{eqnarray*}
The second equivalence happens because $m,m',n,n'$ are all nc analytic 
and each is degree 1 in $h$. \qed 

\section{Complex Hessian as a Sum of Squares} \label{sec:ComplexHessSOS}

Assuming nc plurisubharmonicity means we have a matrix positive nc complex hessian. 
This leads to a sum of squares representation for the nc complex hessian. 

The next lemma follows the proof of Proposition 4.1 in \cite{HM04} with the 
nc hessian now replaced by the nc complex hessian. 

\begin{lemma} \label{lem:deginh}
If $p$ is an nc symmetric plush polynomial then the complex hessian, $q$, of $p$ can be written as
\begin{equation*}
q(x,x^T)[h,h^T] = \sum_{j=1}^{m} r_j^T r_j
\end{equation*}
where each $r_j$ is an nc polynomial that is homogeneous of degree 1 in $h$ (or $h^T$).
\end{lemma}

\proof
Since $p$ is nc plush, that means $q(X,X^T)[H,H^T] \succeq 0$ for every $X$ and $H$. 
By Theorem \ref{thm:Helton}, $q(x,x^T)[h,h^T]$ is a sum of squares. That means we can 
write $q$ as
\begin{equation*}
q(x,x^T)[h,h^T] = \sum_{j=1}^{m} r_j^T r_j
\end{equation*}
where each $r_j$ is a polynomial in $x$, $x^T$, $h$, and $h^T$.
Write
\begin{equation*}
r_j = \sum_{w \in Mon(x,x^T,h,h^T)} {r_j(w) w}
\end{equation*}
where $Mon(x,x^T,h,h^T)$ is the set of monomials in the given variables and where all
but finitely many of the $r_j(w) \in \mathbb{R}$ are 0. Let $deg_h(r)$ denote the degree
of $r$ in $h$ (and $h^T$) and let $deg_x(r)$ denote the degree of $r$ in $x$ (and $x^T$).
Let
\begin{eqnarray*}
d_h &=& \max\{deg_h(r_j) : j\} \\
d_x &=& \max\{deg_x(w) : \exists j \mbox{ s.t. } r_j \mbox{ contains } w \mbox{ and } deg_h(w) = d_h\} \\
S_{d_x,d_h} &=& \{w: r_j \mbox{ contains } w \mbox{ for some } j, deg_h(w) = d_h, deg_x(w) = d_x \}.
\end{eqnarray*}
The portion of $q$ homogeneous of degree $2d_h$ in $h$ and $2d_x$ in $x$ is
\begin{equation*}
\mathcal{Q} = \sum_{\{j=1,\ldots,m, v,w \in S_{d_x,d_h}\}} {r_j(v) r_j(w) v^T w}.
\end{equation*}
Since for $v_j, w_j \in S_{d_x,d_h}$, $v_1^T w_1 = v_2^T w_2$ can occur if and only if $v_1 = v_2$
and $w_1 = w_2$, we see that $Q \ne 0$ and thus $deg_h(q) = 2d_h$. Since $q$ has degree 2 in $h$ and
$h^T$, we obtain $2d_h = 2$ which implies $d_h = 1$.
\qed 

Since we know $q$ is positive, Theorem \ref{thm:Helton} allows us to represent $q$ 
as a sum of squares, $q = \sum r_j^T r_j$. We wish to show that these $r_j$ are either 
analytic or antianalytic.

\begin{theorem} \label{thm:SoSanalyticantianalytic}
If $p$ is an nc symmetric plush polynomial then the complex hessian, $q$, of $p$ can be written as
\begin{equation*}
q(x,x^T)[h,h^T] = \sum_{j=1}^{m} r_j^T r_j
\end{equation*}
where each $r_j$ is either analytic or antianalytic.
\end{theorem}

\proof
Since $p$ is assumed plush we get that $q(X,X^T)[H,H^T] \succeq 0$ for all $X \in (\mathbb{R}^{n \times n})^g$. 
Again, by Theorem \ref{thm:Helton}, we get that $q$ is a sum of squares,
$$ q = \sum {r_j^T r_j}. $$
By Lemma \ref{lem:deginh}, each $r_j$ is homogeneous of degree 1 in $h$ or $h^T$. 
We wish to show that each $r_j$ is either nc analytic or nc antianalytic. Consider all monomials in 
the $r_i$'s of the form
\beq \label{eq:forbidA}
L h_j^T M x_k N
\eeq
or of the form
\beq \label{eq:forbidB}
L x_j^T M h_k N
\eeq
\beq \label{eq:forbidC}
L x_k M h_j^T N
\eeq
or of the form 
\beq \label{eq:forbidD}
L h_k M x_j^T N.
\eeq

Here, $L, M, N$ are monomials in $x$ and $x^T$. The theorem being false is equivalent to some such 
monomial existing and we say these are monomials of the \textit{offending form}. This is easy to check 
just by comparing the form of each offending monomial to $r_j$ being nc analytic or nc antianalytic. 

We now focus on offending monomials of the highest degree (over all offending monomials). \\

\underline{\textbf{Case 1:}} Suppose that the offending monomial of highest degree is of the form 
$L h_j^T M x_k N$. Without loss of generality, say this monomial occurs in $r_1$. Then $r_1^T r_1$ 
contains the monomial
$$ m := N^T x_k^T M^T h_j L^T L h_j^T M x_k N. $$
We claim that this monomial, $m$, appears in $q$. To be cancelled, $\tilde{L} h_j^T M x_k N$ must 
appear in some $r_\ell$ where $\tilde{L}$ factors $L$ or $L$ factors $\tilde{L}$. This implies that either 
$r_\ell^T r_\ell$ contains a monomial of the form $w^T w$, where $w$ is of the offending form and $w$ 
has higher degree than $m$, or we must have $\tilde{L} = L$. The first option would contradict the highest 
degree assumption of $m$ so we must have $\tilde{L} = L$. In this case, the coefficient of $m$ arising from 
$r_\ell^T r_\ell$ is positive so no cancellation occurs.

Observe that $m$ in $q$ contains many Levi-differentially wed monomials in $q$. For example, 
$$ N^T x_k^T M^T x_j L^T L h_j^T M h_k N $$
is contained in $q$, so it appears in some square, say $r_k^T r_k$. Thus, $r_k^T$ contains 
$N^T x_k^T M^T x_j L^T L h_j^T$ (or $N^T x_k^T M^T x_j L^T L h_j^T M$ ) which is of the 
offending form \eqref{eq:forbidC}. But this monomial is longer than the longest offending monomial 
we selected; namely, $m$. This is a contradiction. 

\underline{\textbf{Case 2:}} Suppose that the offending monomial of highest degree is of the form 
$L x_k M h_j^T N$. Without loss of generality, say this monomial occurs in $r_1$. Then $r_1^T r_1$ 
contains the monomial
$$ m := N^T h_j M^T x_k^T L^T L x_k M h_j^T N. $$
We claim that this monomial, $m$, appears in $q$. To be cancelled, $\tilde{K} h_j^T N$ must 
appear in some $r_\ell$ where $\tilde{K}$ factors $L x_k M$ or $L x_k M$ factors $\tilde{K}$. 
This implies that either 
$r_\ell^T r_\ell$ contains a monomial of the form $w^T w$, where $w$ is of the offending form and $w$ 
has higher degree than $m$, or we must have $\tilde{K} = L x_k M$. The first option would contradict the highest 
degree assumption of $m$ so we must have $\tilde{K} = L x_k M$. In this case, the coefficient of $m$ arising from 
$r_\ell^T r_\ell$ is positive so no cancellation occurs. 

Observe that $m$ in $q$ contains many Levi-differentially wed monomials in $q$. For example, 
$$ N^T x_j M^T x_k^T L^T L h_k M h_j^T N $$
is contained in $q$, so it appears in some square, say $r_k^T r_k$. Thus, $r_k^T$ contains 
$N^T x_j M^T x_k^T L^T L h_k$ (or $N^T x_j M^T x_k^T L^T L h_k M$ ) which is of the 
offending form \eqref{eq:forbidB}. But this monomial is longer than the longest offending 
monomial we selected; namely, $m$. This is a contradiction. 

\underline{\textbf{Case 3:}} This case concerns $L h_k M x_j^T N$ and the argument is parallel 
to that in Case 1.  

\underline{\textbf{Case 4:}} This case concerns $L x_j^T M h_k N$ and the argument is parallel 
to that in Case 2.  
\qed 

\section{Proof of Main Results} \label{sec:ProofMainResult}

We now prove our main theorem which we now recall from Section \ref{subsec:MainResult}.

\begin{theorem} \label{thm:pluriRep}
An nc symmetric polynomial $p$ in free variables is nc plush 
if and only if $p$ can be written in the form
\beq \label{eq:pluriRep}
p = \sum f_j^T f_j + \sum k_j k_j^T + F + F^T
\eeq
where the sums are finite and each $f_j$, $k_j$, $F$ is nc analytic.
\end{theorem}

\proof If $p$ has the form given in Equation \eqref{eq:pluriRep} then $q(x,x^T)[h,h^T]$,
the nc complex hessian of $p$, is
\begin{eqnarray*}
q &=& \sum {(f_j^T)_{x^T}[h^T] (f_j)_x[h]} + \sum {(k_j)_x[h] (k_j)_{x^T}[h^T]} \\
&=& \sum {(f_j)_{x}[h]^T (f_j)_x[h]} + \sum {(k_j)_x[h] (k_j)_{x}[h]^T},
\end{eqnarray*}
which is a finite sum of squares. Hence $q(X,X^T)[H,H^T] \succeq 0$ for all $X$ and $H$. 
Thus $p$ is nc plush. 

Suppose $p$ is nc plush. By Theorem \ref{thm:SoSanalyticantianalytic}, we write
$$ 
q = \sum {r_j^T r_j} 
$$
where $r_j$ is analytic or antianalytic homogeneous of degree one in $h$ or $h^T$. 
In view of Theorem \ref{thm:intGvars}, we now show that each $r_j$ is integrable. 
Suppose $m$ is a monomial in $r_j$ 
and that $m'$ is any monomial 1-differentially wed to $m$ (other than $m$). 
We shall now show that $m'$ occurs in $r_j$ with the same coefficient as $m$. 

To do this, suppose $r_j$ contains $C_j m + C_j' m'$ for some $j$'s. Note that
$C_j'$ may certainly be 0. Then $r_j^T r_j$ must contain the terms $C_j^2 m^T m$, $C_j'^2 m'^T m'$,
$C_j C_j' m^T m'$, $C_j' C_j m'^T m$.  By summing over all $j$ such that $r_j$ contains
$C_j m + C_j' m'$ we get that $q$ must have the terms
$$
(\sum_j C_j^2 ) m^T m, \ \
(\sum_j C_j C_j' ) m^T m',  \ \
(\sum_j C_j^{\prime^2} ) m'^T m'
$$
which, by Lemma \ref{lem:monwed}, are Levi-differentially wed. Thus all 3 coefficients
are equal. This means we have
\beq \label{eq:equalcoeffs}
\sum_j C_j^2 = \sum_j C_j C_j' = \sum_j C_j^{\prime^2}.
\eeq
The Cauchy Schwartz inequality gives
\beq \label{eq:CS}
(\sum_j C_j C_j' )^2
\leq
(\sum_j C_j^2 )(\sum_j C_j^{\prime^2} )
\eeq
and Equation \eqref{eq:equalcoeffs} implies we have equality in Equation \eqref{eq:CS}.
This means we have $C_j = \alpha C_j'$ for all $j$. Then we get
$$
\sum_j C_j C_j' = \sum_j \alpha C_j^{\prime^2} = \alpha \sum_j C_j^{\prime^2}
$$
and by Equation \eqref{eq:equalcoeffs}, we get $\alpha = 1$. Hence $C_j = C_j'$ for all
$j$. This means that $r_j$ contains $C_j m$ if and only if it contains $C_j m'$ where
$m$ and $m'$ are any two 1-differentially wed monomials. 

Since $m$ and $m'$ are arbitrary 1-differentially wed monomials,
we get that, by Theorem \ref{thm:intGvars}, $r_j$ is integrable. 
We integrate it to get $f_j$ in Equation \eqref{eq:pluriRep} if $r_j$ 
is nc analytic and $k_j^T$ in Equation \eqref{eq:pluriRep} if $r_j$ 
is nc antianalytic. We note that there are other antiderivatives for the $r_j$ (for example, 
$f_j + x^T$) but when $r_j$ is nc analytic (resp. nc antianalytic) we only care about the 
nc analytic (resp. nc antianalytic) ones.

Define
$$
\tp := \sum f_j^T f_j + \sum k_j k_j^T.
$$
By construction, $\tp$ is a sum of hereditary and antihereditary squares. Also note that
the nc complex hessian of $\tp$ is equal to the nc complex hessian of $p$. Apply Lemma 
\ref{lem:leviUnique} to finish the proof. \qed 

Now we prove the uniqueness of the representation of an nc symmetric 
plush polynomial. We recall Theorem \ref{thm:MainThmUnique} from Section 
\ref{subsec:MainResult}:

\begin{theorem} \label{thm:Uniqueness}
Let $p$ be an nc symmetric polynomial in free variables that is nc plurisubharmonic 
and let
\begin{eqnarray*}
\tN &:=& \min \{N : p = \sum_{j=1}^{N}{f_j^T f_j} + \sum_{j=1}^{M}{k_j k_j^T} + F + F^T\} \\
\tM &:=& \min \{M: p = \sum_{j=1}^{N}{f_j^T f_j} + \sum_{j=1}^{M}{k_j k_j^T} + F + F^T\}.
\end{eqnarray*}
Then, we can represent $p$ as 
$$
p = \sum_{j=1}^{\tN}{\tf_j^T \tf_j} + \sum_{j=1}^{\tM}{\tk_j \tk_j^T} + \tilde{F} +\tilde{F}^T
$$
and if $N$ and $M$ are integers such that $N \ge \tN$, $M \ge \tM$ and 
$$
p = \sum_{j=1}^{N}{f_j^T f_j} + \sum_{j=1}^{M}{k_j k_j^T} + F + F^T \ ,
$$
then there exist isometries $U_1 : \mathbb{R}^{\tN} \longrightarrow \mathbb{R}^{N}$ and 
$U_2 : \mathbb{R}^{\tM} \longrightarrow \mathbb{R}^{M}$ such that 
$$
\left( \begin{array}{c}
f_1 \\
\vdots \\
f_N
\end{array} \right)
=
U_1 
\left( \begin{array}{c}
\tf_1 \\
\vdots \\
\tf_{\tN}
\end{array} \right)
+ 
\vec{c}_1
\ \ \ \mbox{and} \ \ \ 
\left( \begin{array}{c}
k_1^T \\
\vdots \\
k_M^T
\end{array} \right)
=
U_2
\left( \begin{array}{c}
\tk_1^T \\
\vdots \\
\tk_{\tM}^T
\end{array} \right)
+ 
\vec{c}_2
$$
where $\vec{c}_1 \in \mathbb{R}^N$ and $\vec{c}_2 \in \mathbb{R}^M$.
\end{theorem}

\proof Suppose $N$ and $M$ are integers such that $N \ge \tN, M \ge \tM$ 
where $p$ can be written as
\beq \label{eq:tempPsh2}
p = \sum_{j=1}^{N}{f_j^T f_j} + \sum_{j=1}^{M}{k_j k_j^T} + F + F^T
\eeq
and suppose $\widehat{M} \ge \tM$ is such that we have  
\beq \label{eq:tempPsh}
p = \sum_{j=1}^{\tN} {\tf_j^T \tf_j} + \sum_{j=1}^{\widehat{M}} {\hat{k}_j \hat{k}_j^T} + F + F^T.
\eeq

Then, the nc complex hessian, $q$, of $p$ based on the representations 
in Equations \eqref{eq:tempPsh} and \eqref{eq:tempPsh2} is 
\beq \label{eq:temp}
q = \sum_{j=1}^{\tN} {(\tf_j)_x[h]^T (\tf_j)_x[h]} + \sum_{j=1}^{\widehat{M}} {(\hat{k}_j)_x[h] (\hat{k}_j)_x[h]^T} 
\eeq
$$
= \sum_{j=1}^{N} {(f_j)_x[h]^T (f_j)_x[h]} + \sum_{j=1}^{M} {(k_j)_x[h] (k_j)_x[h]^T}.
$$
We define $\purehered{q}$ as the \textbf{purely hereditary} part of $q$ to be 
all of the terms that contain $h^T$ to the left of $h$ and we 
define $\pureantihered{q}$ as the \textbf{purely antihereditary} part of $q$ to be 
all of the terms that contain $h$ to the left of $h^T$. 

First, from Equation \eqref{eq:temp}, consider the purely hereditary part of the nc complex hessian, 
$$
\purehered{q} = \sum_{j=1}^{\tN} {(\tf_j)_x[h]^T (\tf_j)_x[h]} = \sum_{j=1}^{N} {(f_j)_x[h]^T (f_j)_x[h]}.
$$
Since $\purehered{q}$ is a sum of squares, it is matrix positive. Hence, 
the Gram representations\footnote{See \cite{PW} as a reference for the 
Gram representation.} for this purely hereditary part of $q$ contain 
$\ell \times \ell$ unique positive semidefinite matrices, $\widetilde{G}$ and 
$G$, that are both of rank $\tN$ such that 
$$
\sum_{j=1}^{\tN} {(\tf_j)_x[h]^T (\tf_j)_x[h]} = y^T \widetilde{G} y 
\ \ \ \mbox{and} \ \ \
\sum_{j=1}^{N} {(f_j)_x[h]^T (f_j)_x[h]} = y^T G y \ ,
$$
where $y$ is an $\ell \times 1$ vector of monomials in $x$ and $h$. 
The purely hereditary nature of $\purehered{q}$ forces $\widetilde{G}$ 
and $G$ to be unique (so, in fact, $\widetilde{G}=G$). 

Since $\widetilde{G}$ is positive semidefinite, we can write $\widetilde{G}$ as 
$\widetilde{G} = \widetilde{W}^T \widetilde{W}$,  
where $\widetilde{W} : \mathbb{R}^{\ell} \longrightarrow \mathbb{R}^{\tN}$ 
is an $\tN \times \ell$ matrix with $\rank(\widetilde{W})=\tN$ such that 
$$
\widetilde{W} y = 
\left( \begin{array}{c}
(\tf_1)_x[h] \\
\vdots \\
(\tf_{\tN})_x[h]
\end{array} \right).
$$

Similarly, we can write $G$ as $G = W^T W$, where $W: \mathbb{R}^{\ell} 
\longrightarrow \mathbb{R}^{N}$ is an $N \times \ell$ matrix with 
$\rank(W) = \tN$ such that 
$$
W y = 
\left( \begin{array}{c}
(f_1)_x[h] \\
\vdots \\
(f_{N})_x[h]
\end{array} \right).
$$
Note that the range of $W$ is an $\tN$-dimensional subspace 
sitting inside of $\mathbb{R}^N$.

Let $\cR^{\tN}$ denote the subspace of $\mathbb{R}^N$ spanned 
by the first $\tN$ coordinates of $\mathbb{R}^N$; i.e., 
$$
\cR^{\tN} = \span \{ e_1, e_2, \ldots, e_{\tN} \}
$$
where $e_i$ is the $i^{th}$ standard basis vector in $\mathbb{R}^N$. 
Then define the $N \times \tN$ matrix $E : \mathbb{R}^{\tN} \longrightarrow 
\mathbb{R}^{N}$ as 
$
E = 
\left( \begin{array}{c}
I_{\tN} \\
0
\end{array} \right)
$
so that
$
E \widetilde{W} y = 
\left( \begin{array}{c}
\widetilde{W} y \\
0
\end{array} \right) 
\in \mathbb{R}^N.
$
We note that if $N = \tN$, then $E = I_{\tN}$.

Let $V : \mathbb{R}^N \longrightarrow \mathbb{R}^N$ be an 
$N \times N$ unitary matrix that maps $\cR^{\tN}$ onto the 
range of $W$ such that 
$$
V \left( \begin{array}{c}
\widetilde{W} y \\
0
\end{array} \right) 
= 
W y.
$$

This implies that $W y = V E \widetilde{W} y$ and note that the 
$N \times \tN$ matrix $U_1 = V E : \mathbb{R}^{\tN} \longrightarrow 
\mathbb{R}^N$ is an isometry and that 
\beq \label{eq:isomInt}
\left( \begin{array}{c}
(f_1)_x[h] \\
\vdots \\
(f_N)_x[h]
\end{array} \right)
=
U_1 
\left( \begin{array}{c}
(\tf_1)_x[h] \\
\vdots \\
(\tf_{\tN})_x[h]
\end{array} \right).
\eeq
Now, we perform nc integration to each nc polynomial in the vectors on 
both sides of Equation \eqref{eq:isomInt}. We do this according to Corollary 
\ref{cor:Constant} to get 
$$
\left( \begin{array}{c}
f_1 \\
\vdots \\
f_N
\end{array} \right) 
= 
U_1 
\left( \begin{array}{c}
\tf_1 \\
\vdots \\
\tf_{\tN}
\end{array} \right) 
+ 
\vec{c}_1
$$
where $\vec{c}_1 \in \mathbb{R}^{N}$.

Similarly, if, at the start of the proof, we assumed $\widehat{N} \ge \tN$ 
is such that we have 
$$
p = \sum_{j=1}^{\widehat{N}}{\hat{f}_j^T \hat{f}_j} + \sum_{j=1}^{\tM}{\tk_j \tk_j^T} + F + F^T \ ,
$$
then, we would have constructed an isometry $U_2 : \mathbb{R}^{\tM} \longrightarrow 
\mathbb{R}^{M}$ such that 
$$
\left( \begin{array}{c}
k_1 \\
\vdots \\
k_M
\end{array} \right)
=
U_2
\left( \begin{array}{c}
\tk_1 \\
\vdots \\
\tk_{\tM}
\end{array} \right)
+
\vec{c}_2
$$ 
where $\vec{c}_2 \in \mathbb{R}^{M}$ and $M$ is as in Equation \eqref{eq:tempPsh2}.

By defining the isometry $U = U_1 \oplus U_2$, we can then write $p$ with the 
minimal number of hereditary and antihereditary squares as 
$$
p = \sum_{j=1}^{\tN}{\tf_j^T \tf_j} + \sum_{j=1}^{\tM}{\tk_j \tk_j^T} + \tilde{F} + \tilde{F}^T.
$$
\qed

\end{document}